\documentclass[11pt]{article}

\usepackage{latexsym}
\usepackage{amssymb}

\newtheorem{thm}{Theorem}

\newtheorem{rmr}{Remark}

\begin{document}
{
\begin{center}
{\Large\bf
The Nevanlinna-type parametrization for the operator Hamburger moment problem.}
\end{center}
\begin{center}
{\bf S.M. Zagorodnyuk}
\end{center}

\section{Introduction.}

Let $\mathcal{H}$ be an arbitrary (not necessarily separable) Hilbert space. 
\textit{The operator Hamburger moment problem} consists of
finding a non-decreasing $[\mathcal{H}]$-valued function $F(t)$, $t\in\mathbb{R}$, $F(0)=0$,
which is strongly left-continuous on $\mathbb{R}$ and
such that
\begin{equation}
\label{f1_1}
\int_\mathbb{R} t^n dF(t) := \mathrm{s.-} \lim\limits_{ {a\to -\infty}\atop{b\to +\infty}  } \int_a^b t^n dF(t) =
S_n,\qquad n\in \mathbb{Z}_+,
\end{equation}
where $\{ S_n \}_{n=0}^\infty$ is a prescribed sequence of bounded operators on $\mathcal{H}$ (called \textit{moments}).
The operator Stieltjes integral $\int_a^b t^n dF(t)$ is understood as a limit of the corresponding Stieltjes-type operator
integral sums in the strong operator topology.
%Namely, let $a=t_0 < t_1 < ... < t_N=b$ be a partition of $[a,b]$, $d = \max\limits_{ 0\leq j\leq N-1} |t_{j+1} - t_j|$ be its diameter,
%and $t_k^*\in [t_k,t_{k+1})$, $0\leq k\leq N-1$ be arbitrary points. Then
%$\int_a^b t^n dF(t) h = \lim_{d\to 0} \sum_{k=0}^{N-1} (t_k^*)^n (F(t_{k+1})-F(t_k)) h$, $h\in\mathcal{H}$.

If the  moment problem~(\ref{f1_1}) has a solution, then $\{ S_n \}_{n=0}^\infty$ is said to be \textit{a moment sequence}.
The  moment problem~(\ref{f1_1}) is said to be {\it determinate} if it has a unique solution and {\it indeterminate}
if it has more than one solution. 
The following conditions:
\begin{equation}
\label{f1_2}
\sum_{k,l=0}^r (S_{k+l} h_k, h_l)_{\mathcal{H}} \geq 0,\quad 
h_j\in\mathcal{H}\ (0\leq j\leq r),\quad
\forall r\in\mathbb{Z}_+, 
\end{equation}
are necessary and sufficient for the solvability of the moment problem~(\ref{f1_1}). 
Observe that the conditions of the solvability 
appeared in~\cite{cit_510_I} without a proof, while in~\cite{citNew_20_K} they appeared with a reference to~\cite{cit_600_Akh} (however we
do not remember such a statement in~\cite{cit_600_Akh}), see also~\cite{citNew_25_LN}.

In the scalar case ($\mathcal{H}=\mathbb{C}$), a description of all solutions of the moment problem~(\ref{f1_1}) can be found, e.~g.,
in~the classical books \cite{cit_600_Akh},~\cite{cit_700_Ber} (for the nondegenerate case) and in~\cite{cit_800_AK} (for
the degenerate case).

The matrix case ($\mathcal{H}=\mathbb{C}_{N\times N}$), introduced by Krein in~\cite{citNew_30_K}, nowadays still contains unanswered questions and 
attracts many researchers,
see recent works~\cite{citNew_50_CH}, \cite{cit_400_L}, \cite{cit_970_B}, \cite{cit_850_Z}, \cite{citNew_100_Z} and historical references therein.

Probably, the first operator moment problems, namely the Stieltjes and Hausdorff operator moment problems, were introduced in~1947 by Krein
and Krasnoselskii in~\cite[Examples 8.1 and 10.1]{citNew_300_KK}.
On page~97 in~\cite{citNew_300_KK} they also noticed that the power moment problem on $(-\infty,\infty)$ admitted the corresponding extension
for the case of the operator moment problem. Krein and Krasnoselskii obtained conditions of the solvability for
the Stieltjes and Hausdorff operator moment problems.
In~1953 Sz.-Nagy, considering another problem, solved an operator power moment problem on a compact subset of the reals, 
see~\cite[pp. 286-288]{citNew_500_SN}.
In~1962 MacNerney characterized moment sequences for the operator Hamburger moment problem \cite{citNew_700_MN}.
Berezansky in 1965 proved the convergence of the series
from the polynomials of the first kind for the operator moment problem~\cite[Ch.7, Section 2]{cit_700_Ber}.
In~1966 Leviatan considered a modified operator Hausdorff moment problem and obtained conditions for its solvability~\cite{citNew_900_L}.
Set
\begin{equation}
\label{f1_3}
\Gamma_n = \left(
\begin{array}{cccc} S_0 & S_1 & \ldots & S_n\\
S_1 & S_2 & \ldots & S_{n+1}\\
\vdots & \vdots & \ddots & \vdots\\
S_n & S_{n+1} & \ldots & S_{2n}\end{array}
\right),\qquad n\in \mathbb{Z}_+.
\end{equation}
Here $\Gamma_n$ may be viewed as an operator on the direct sum of $n+1$ copies of $\mathcal{H}$.
Under a condition of the strict positivity of $\Gamma_n$, $\forall n\in \mathbb{Z}_+$,
all solutions of the operator Hamburger moment problem~(\ref{f1_1}) were parametrized by
Ilmushkin in~\cite{cit_510_I}, and, in a different way, by 
Ilmushkin and Aleksandrov in~\cite{cit_520_IA}. 
However we should notice that the
formula, connecting spectral functions of the corresponding symmetric operator and solutions
of the moment problem, was stated without proof in~\cite[formula (2.3)]{cit_520_IA} and referred
papers~\cite[Theorem 2]{cit_510_I},\cite[Theorem 1]{cit_530_I}. 
Moreover, it was stated without proof
that the latter formula holds true for the general case~\cite[pp. 77-78]{cit_510_I}.

\noindent
Probably a first description of all solutions of the operator Hamburger moment problem, without additional assumptions,
was given by Kheifets~\cite{citNew_20_K}. He used the well-known abstract interpolation problem (AIP),
introduced by Katsnelson, Kheifets and Yuditskii. This description was not given by a Nevanlinna-type
formula and it had a more complicated structure~\cite{citNew_20_K}:
\begin{equation}
\label{f1_4_1}
\int_{-\infty}^\infty \frac{\sigma(dx)}{x-z} =
i\frac{1+w(\zeta)}{1-w(\zeta)},\qquad z = i\frac{1+\zeta}{1-\zeta},
\end{equation}
\begin{equation}
\label{f1_4_2}
w = s_0 + s_2 (1_{N_2} - \omega s)^{-1} \omega s_1,
\end{equation}
where $\sigma(dx)$ is a solution of the moment problem, $S = \left(\begin{array}{cc} s & s_1\\
s_2 & s_0\end{array}\right)$ is the scattering function of the AIP, and $\omega(\zeta)$ is an arbitrary
analytic operator-valued contractive function on $\mathbb{D}$, $\omega(\zeta):N_1\rightarrow N_2$
($N_1$,$N_2$ be some Hilbert spaces). 
We do not see an easy way to simplify the description~(\ref{f1_4_1})-(\ref{f1_4_2}) to get
a Nevanlinna-type parametrization.

Here we focused our attension mainly on the history of the (full) operator Hamburger moment problem. The truncated operator moment
problems and various modifications were also intensively studied, see, e.~g., \cite{citNew_1000_A},
\cite{cit_520_IA}, \cite{citNew_1400_A}, \cite{citNew_1700_Z} and references therein.

Our aim here is to obtain a Nevanlinna-type parametrization for the operator Hamburger moment problem~(\ref{f1_2}).
This description will be given by a linear fractional transformation with bounded operators as coefficients.
No additional assumptions besides solvability of the moment problem will be posed.

We shall apply the operator approach to the moment problem~(\ref{f1_1}).
As far as we know, in 1940 Naimark laid foundation of the operator approach to moment problems in his paper~\cite{citNew_12000_Neumark}.
Namely, he described all solutions of the Hamburger moment problem
in terms of spectral functions of the operator defined by the Jacobi matrix~\cite[pp. 303-305]{citNew_12000_Neumark}.
Three years later, in 1943 Neumark, using his description of the generalized resolvents, derived
Nevanlinna's formula for all solutions of the Hamburger moment problem~\cite[pp. 292-294]{citNew_14000_Neumark}.
In 1947 Krein and Krasnoselskii presented their version of the operator approach to the Hamburger moment
problem in~\cite{citNew_300_KK}. An abstract operator approach to
the Nevanlinna-Pick problem was proposed by
Sz\"okefalvi-Nagy and Koranyi in~\cite{citNew_4280_SK}, \cite{citNew_4290_SK}.
Our approach to the operator moment problem~(\ref{f1_1}) is close to the ideas of
Krein, Krasnoselskii, Sz\"okefalvi-Nagy and Koranyi.
Similar ideas were applied in our paper~\cite{citNew_1700_Z} for the case of the operator trigonometric moment
problem. However, in the case of the moment problem on the real line there appear new difficulties.

The first problem, which we shall meet in the next section, is to define a version of the $L_2$ space for
$\mathcal{H}$-valued functions, equipped with enough functions for solving the moment problem~(\ref{f1_1}).
For the case of a \textit{separable} Hilbert space $H$ and an operator-valued function $T(\lambda)$ 
(subject to some conditions) on $\mathbb{R}$, the space 
$L_2(H;(-\infty,\infty),dT(\lambda))$ was introduced in~\cite[p. 556]{cit_700_Ber}. Let us recall this definition.
Denote by $C_{00}(H;(-\infty,\infty))$ the set of all strongly continuous finitely supported  $H$-valued functions $f(\lambda)$
($-\infty < \lambda <\infty$) which values are situated in a finite-dimensional subspace (depending on $f$) of $H$.
For $f,g\in C_{00}(H;(-\infty,\infty))$ the scalar product is given by the following formula:
\begin{equation}
\label{f1_5}
(f,g) :=  \int_{-\infty}^\infty (dT(\lambda) f(\lambda), g(\lambda)),
\end{equation}
where, according to~\cite{cit_700_Ber}, the integral is a weak limit of the corresponding Riemann-Stieltjes integral sums for
finite partitions of the real line.
%The first question which appears: how the finite partitions of the (infinite) real line become arbitrarily fine?
%This is not answered in~\cite{cit_700_Ber}. 
%Also important for our purposes is another question. 
The space $L_2(H;(-\infty,\infty),dT(\lambda))$ is defined by standard procedures of the factorization and the completion from
$C_{00}(H;(-\infty,\infty))$. It is noticed in~\cite[p. 557]{cit_700_Ber} that if the function $f(\lambda)$ is 
strongly continuous on $(-\infty,\infty)$ and
$\int_{-\infty}^\infty (dT(\lambda) f(\lambda), f(\lambda)) < \infty$, 
then $f\in L_2(H;(-\infty,\infty),dT(\lambda))$.
What was meant by this? Probably, it was understood that for such $f(\lambda)$ there exists a fundamental
sequence (of classes of the equivalence) of functions $f_n(\lambda)$ from $C_{00}(H;(-\infty,\infty))$ ($n\in\mathbb{N}$),
which tends to $f(\lambda)$ as $n\rightarrow\infty$: $(f-f_n,f-f_n)\rightarrow 0$. 
Observe that $f_n$ should be constructed finitely supported, strongly continuous and with values in a finite-dimensional subspace.
No discussion on this construction can be found in~\cite{cit_700_Ber}.
Moreover, since it should use a finite-dimensional subspace approximation of $f(\lambda)$, it seems to be not
applicable for the case of a non-separable $\mathcal{H}$.

The above abstract identification of functions with fundamental sequences seems to be not transparent and not convenient
in applications.
We prefer to deal with usual functions (or classes of equivalent functions) as elements of $L_2$ type space.
Thus, in Section~2 we shall define a space $L_2(\mathcal{H},dF(t))$ of $\mathcal{H}$-valued functions on $\mathbb{R}$ with enough set of
basic functions for our purposes. The operator $\mathcal{A}$ of the multiplication by an independent variable in $L_2(\mathcal{H},dF(t))$
is studied, as well.

In Section~3 we shall study the moment problem~(\ref{f1_1}) introducing an abstract Hilbert space $H$ and a symmetric operator $A$
in it (by means of given moments). The spectral functions of $A$ generate solutions of the moment problem. 
The space $L_2(\mathcal{H},dF(t))$ and the operator $\mathcal{A}$ play an important role in establishing that each solution
of the moment problem is generated by a spectral function of $A$.
In passing to Nevanlinna-type description of solutions important tools are
Chumakin's formula (\cite{cit_950_Ch}) for the generalized resolvents of a closed isometric operator (applied for Cayley's transformation of $A$) 
and Frobenius-type formula for operator matrices. 
As differ from the case of the truncated operator trigonometric moment problem, we need to
shift the initial vector $x_0\in H$ to get into a block of the Frobenius-type formula, similar to the matrix case in~\cite{citNew_100_Z}.
After some computations we shall come to a desired Nevanlinna-type parametrization
of all solutions of the operator Hamburger moment problem.
The coefficients of the corresponding operator linear fractional transformation are bounded operators.
%Conditions for the determinacy of the moment problem (in terms of the defect numbers of $A$) are obtained,
%as well.

%

\noindent
{\bf Notations.}  
As usual, we denote by $\mathbb{R}, \mathbb{C}, \mathbb{N}, \mathbb{Z}, \mathbb{Z}_+$,
the sets of real numbers, complex numbers, positive integers, integers and non-negative integers,
respectively; $\mathbb{D} = \{ z\in \mathbb{C}:\ |z|<1 \}$, 
$\mathbb{R}_e = \mathbb{C}\backslash\mathbb{R}$,
$\mathbb{C}_+ = \{ z\in \mathbb{C}:\ \mathop{\rm Im}\nolimits z > 0\}$.
By $\mathfrak{B}(\mathbb{R})$ we mean the set of all Borel subsets of the real line.

\noindent
In this paper Hilbert spaces are not necessarily separable, operators in them are
supposed to be linear.

\noindent
If H is a Hilbert space then $(\cdot,\cdot)_H$ and $\| \cdot \|_H$ mean
the scalar product and the norm in $H$, respectively.
Indices may be omitted in obvious cases.
For a linear operator $A$ in $H$, we denote by $D(A)$
its  domain, by $R(A)$ its range, and $A^*$ means the adjoint operator
if it exists. If $A$ is invertible then $A^{-1}$ means its
inverse. $\overline{A}$ means the closure of the operator, if the
operator is closable. If $A$ is bounded then $\| A \|$ denotes its
norm.
For a set $M\subseteq H$
we denote by $\overline{M}$ the closure of $M$ in the norm of $H$.
For an arbitrary set of elements $\{ x_n \}_{n\in I}$ in
$H$, we denote by $\mathop{\rm Lin}\nolimits\{ x_n \}_{n\in I}$
the set of all linear combinations of elements $x_n$,
and $\mathop{\rm span}\nolimits\{ x_n \}_{n\in I}
:= \overline{ \mathop{\rm Lin}\nolimits\{ x_n \}_{n\in I} }$.
Here $I$ is an arbitrary set of indices.
By $E_H$ we denote the identity operator in $H$, i.e. $E_H x = x$,
$x\in H$. In obvious cases we may omit the index $H$. If $H_1$ is a subspace of $H$, then $P_{H_1} =
P_{H_1}^{H}$ is an operator of the orthogonal projection on $H_1$
in $H$.
By $[H]$ we denote the set of all bounded operators on $H$.
%For a closed isometric operator $V$ in $H$ we denote:
%$M_\zeta(V) = (E_H - \zeta V) D(V)$, $N_\zeta(V) = H\ominus M_\zeta(V)$, $\zeta\in \mathbb{T}_e$; $M_\infty(V)=R(V)$, $N_\infty(V)= H\ominus R(V)$.
For a closed symmetric operator $B$ in $H$ we denote:
$\mathcal{M}_z(B) = (B - z E_H) D(B)$, $\mathcal{N}_z(B) = H\ominus \mathcal{M}_z(B)$, $z\in \mathbb{R}_e$.
If $B$ is self-adjoint, then we set $R_z(B) = (B-z E_H)^{-1}$, $z\in\mathbb{R}_e$.

\noindent
By $\mathcal{S}(D;N,N')$ we denote a class of all analytic in a domain $D\subseteq \mathbb{C}$
operator-valued functions $F(z)$, which values are linear non-expanding operators mapping the whole
$N$ into $N'$, where $N$ and $N'$ are some Hilbert spaces.

\section{The space $L_2(\mathcal{H},dF(t))$.}

Let $\mathcal{H}$ be an arbitrary, not necessarily separable, Hilbert space. 
Let $F(t)$ be a non-decreasing $[\mathcal{H}]$-valued function, $t\in\mathbb{R}$, $F(0)=0$,
which is strongly left-continuous on $\mathbb{R}$.
%and
%such that there exists
%\begin{equation}
%\label{f2_3}
%\mathrm{s.-} \lim\limits_{ {a\to -\infty}\atop{b\to +\infty}  } (F(b)-F(a)), 
%\end{equation}
%and the following operator $F_0$:
%\begin{equation}
%\label{f2_5}
%F_0 := F(+\infty) - F(-\infty) := \mathrm{s.-} \lim\limits_{ {a\to -\infty}\atop{b\to +\infty}  } (F(b)-F(a))
%\end{equation}
%is a bounded operator on $\mathcal{H}$.
The space $\mathcal{H}$ and the function $F(t)$ will be fixed throughout this section.

\begin{rmr}
\label{r2_1}
In this section we shall not consider the moment problem~(\ref{f1_1}). Thus, the space $\mathcal{H}$ and the function $F(t)$
are not related to any moment problem. We think that our further constructions have some
interest themselves. Of course, our main aim here is to apply these constructions to the moment problem in the next section.
\end{rmr}

Denote by $\Phi = \Phi(\mathcal{H},dF(t))$ a set of all $\mathbb{C}$-valued continuous functions $\varphi(t)$, $t\in\mathbb{R}$, such that
\begin{equation}
\label{f2_7}
\int_{\mathbb{R}} |\varphi(t)|^2 d (F(t) h,h)_{\mathcal{H}} 
:=
\lim\limits_{ {a\to -\infty}\atop{b\to +\infty}  } \int_a^b |\varphi(t)|^2 d (F(t) h,h)_{\mathcal{H}}  
< \infty,\qquad \forall h\in\mathcal{H}.
\end{equation}
%In particular, we see that all continuous, bounded $\mathbb{C}$-valued functions belong to $\Phi(\mathcal{H},dF(t))$.

Denote by $\Theta = \Theta(\mathcal{H},dF(t))$ a set of all $\mathcal{H}$-valued functions $f(t)$, $t\in\mathbb{R}$, which
admit a representation of the following form:
\begin{equation}
\label{f2_9}
f(t) = \sum_{k=0}^r \varphi_k(t) h_k,\qquad  \varphi_k\in\Phi,\ h_k\in\mathcal{H},\ r\in\mathbb{Z}_+.
\end{equation}

For arbitrary two functions $f(t)$ and $g(t)$ from $\Theta(\mathcal{H},dF(t))$ we set
\begin{equation}
\label{f2_11}
\Psi(f,g) := \int_{\mathbb{R}} (dF(t) f(t), g(t))_{\mathcal{H}} :=  
\lim\limits_{ {a\to -\infty}\atop{b\to +\infty}  } \int_a^b (dF(t) f(t), g(t))_{\mathcal{H}},
\end{equation}
where
\begin{equation}
\label{f2_15}
\int_a^b (dF(t) f(t), g(t))_{\mathcal{H}} := \lim_{\delta\rightarrow +0} \sum_{j=0}^{N-1} (
( F(t_{j+1}) - F(t_j) ) f(t_j^*), g(t_j^*)
)_{\mathcal{H}}.
\end{equation}
Here $\delta$ is the diameter of a partition of $[a,b]$:
\begin{equation}
\label{f2_16}
a=t_0 < t_1 < t_2 < ... < t_{N}=b, 
\end{equation}
and $t_j^*\in [t_j, t_{j+1})$ are arbitrary points.
As usual, it is understood that the limit in~(\ref{f2_15}) does not depend on the choice of partitions and points $t_j^*$.

\noindent
Let us check that limits in~(\ref{f2_15}) and~(\ref{f2_11}) do exist for all functions $f,g\in \Theta(\mathcal{H},dF(t))$,
and calculate the value of $\Psi(f,g)$. 
Let $f(t)$ have the form~(\ref{f2_9}) and
\begin{equation}
\label{f2_17}
g(t) = \sum_{l=0}^q \psi_l(t) g_l,\qquad  \psi_l\in\Phi,\ g_l\in\mathcal{H},\ q\in\mathbb{Z}_+.
\end{equation}
By the substitution of representations~(\ref{f2_9}) and~(\ref{f2_17}) into expressions in~(\ref{f2_15}) we get
$$ \int_a^b (dF(t) f(t), g(t))_{\mathcal{H}} = $$
$$ =
\lim_{\delta\rightarrow +0}
\sum_{k=0}^r \sum_{l=0}^q
\sum_{j=0}^{N-1} \varphi_k(t_j^*) \overline{ \psi_l(t_j^*) } 
\left(
(F(t_{j+1})h_k, g_l)_{\mathcal{H}} - (F(t_{j})h_k, g_l)_{\mathcal{H}}
\right)
= $$
\begin{equation}
\label{f2_19}
= \sum_{k=0}^r \sum_{l=0}^q \int_a^b \varphi_k(t) \overline{ \psi_l(t) } d (F(t)h_k, g_l)_{\mathcal{H}}.
\end{equation}
Then
$$ \Psi(f,g) =  
\lim\limits_{ {a\to -\infty}\atop{b\to +\infty}  } 
\sum_{k=0}^r \sum_{l=0}^q \int_a^b \varphi_k(t) \overline{ \psi_l(t) } d (F(t)h_k, g_l)_{\mathcal{H}} = $$
\begin{equation}
\label{f2_25}
= \sum_{k=0}^r \sum_{l=0}^q \int_{\mathbb{R}} \varphi_k(t) \overline{ \psi_l(t) } d (F(t)h_k, g_l)_{\mathcal{H}},
\end{equation}
where
\begin{equation}
\label{f2_27}
\int_{\mathbb{R}} \varphi_k(t) \overline{ \psi_l(t) } d (F(t)h_k, g_l)_{\mathcal{H}}
:=
\lim\limits_{ {a\to -\infty}\atop{b\to +\infty}  } 
\int_a^b \varphi_k(t) \overline{ \psi_l(t) } d (F(t)h_k, g_l)_{\mathcal{H}}.
\end{equation}
It remains to explain why the latter limit exists.
In fact, by the polarization formula we may write:
\begin{equation}
\label{f2_28}
\int_a^b \varphi_k(t) \overline{ \psi_l(t) } d (F(t)h_k, g_l)_{\mathcal{H}} =
\sum_{j=1}^4 c_j \int_a^b \varphi_k(t) \overline{ \psi_l(t) } d (F(t) u_j, u_j)_{\mathcal{H}},
\end{equation}
where $ c_{1,2} = \pm \frac{1}{4}, c_{3,4} = \pm \frac{i}{4}$,
$u_{1,2} = h_k \pm g_l, u_{3,4} = h_k \pm i g_l$.
By condition~(\ref{f2_7}) we conclude that there exist the following limits:
$$ \lim\limits_{ {a\to -\infty}\atop{b\to +\infty}  } 
\int_a^b \varphi_k(t) \overline{ \psi_l(t) } d (F(t) u_j, u_j)_{\mathcal{H}}, $$ 
and therefore the required limit exists.

It is clear that $\Theta(\mathcal{H},dF(t))$ is a complex linear vector space with usual operations of the sum and the multiplication
by a complex scalar.
By~(\ref{f2_11}), (\ref{f2_15}) we see that the functional $\Psi(f,g)$ is sesquilinear,
$\overline{ \Psi(f,g) } = \Psi(g,f)$,
and $\Psi(f,f)\geq 0$ ($f,g\in\Theta$).
As usual, we put $f$ and $g$ from $\Theta$ to the same class of the equivalence, if
$\Psi(f-g,f-g) = 0$. 
The set of classes of the equivalence we denote by $\widetilde\Theta = \widetilde\Theta(\mathcal{H},dF(t))$.
By the completion of $\widetilde\Theta$ we obtain a Hilbert space, which we denote by $L_2(\mathcal{H},dF(t))$.
In what follows we denote by $[f]$ the class of the equivalence in $L_2(\mathcal{H},dF(t))$ which contains $f\in\Theta$.

Denote by $\Phi_0 = \Phi_0(\mathcal{H},dF(t))$  a set of all continuous $\mathbb{C}$-valued functions $\varphi(t)$
such that $\varphi(t) = 0$ for all $t$ lying outside some finite interval $(c,d)$, where $(c,d)$ may depend on $\varphi$
(i.e. a set of all continuous finitely supported $\mathbb{C}$-valued functions).
Of course, we have inclusion $\Phi_0(\mathcal{H},dF(t))\subseteq\Phi(\mathcal{H},dF(t))$.
For an arbitrary function $\varphi(t)\in \Phi(\mathcal{H},dF(t))$ 
there exist functions $\varphi^{[n]}(t)\in\Phi_0(\mathcal{H},dF(t))$, $n\in\mathbb{N}$, such that for
an arbitrary $h\in\mathcal{H}$ we have:
$$ \int_{\mathbb{R}} |\varphi_(t) - \varphi^{[n]}(t) |^2 d(F(t)h,h)_{\mathcal{H}} =
$$
\begin{equation}
\label{f2_29}
= \lim\limits_{ {a\to -\infty}\atop{b\to +\infty}  } \int_a^b 
|\varphi_(t) - \varphi^{[n]}(t) |^2 d(F(t)h,h)_{\mathcal{H}} 
\rightarrow 0,\quad \mbox{as } n\rightarrow\infty.
\end{equation}
For example, we may choose
\begin{equation}
\label{f2_31}
\varphi^{[n]}(t) =
\left\{
\begin{array}{cccc} \varphi(t), & t\in [-n,n]\\
\varphi(t) (-t+n+1), & t\in [n,n+1]\\
\varphi(t) (t+n+1), & t\in [-n-1,-n]\\
0, & t\in (-\infty,-n-1]\cup[N+1,\infty)\end{array}
\right..
\end{equation}

Denote by $\Theta_0 = \Theta_0(\mathcal{H},dF(t))$ a set of all $\mathcal{H}$-valued functions $f(t)$, $t\in\mathbb{R}$, which
admit a representation of the following form:
\begin{equation}
\label{f2_33}
f(t) = \sum_{k=0}^r \varphi_k(t) h_k,\qquad  \varphi_k\in\Phi_0,\ h_k\in\mathcal{H},\ r\in\mathbb{Z}_+.
\end{equation}
Notice that $\Theta_0\subseteq\Theta$. Set
$$ \widetilde\Theta_0 = \widetilde\Theta_0(\mathcal{H},dF(t)) = 
\left\{
[f(t)],\ f(t)\in \Theta_0(\mathcal{H},dF(t))
\right\}, $$
where $[f]$ means the class of the equivalence in $L_2(\mathcal{H},dF(t))$ which contains $f$.
It is clear that $\widetilde\Theta_0 \subseteq \widetilde\Theta$.
Let us check that the set $\widetilde\Theta_0(\mathcal{H},dF(t))$ is dense in $L_2(\mathcal{H},dF(t))$.
Since $\widetilde\Theta(\mathcal{H},dF(t))$ is dense in $L_2(\mathcal{H},dF(t))$ by the construction of $L_2(\mathcal{H},dF(t))$, it is enough
to approximate an arbitrary element from $\widetilde\Theta(\mathcal{H},dF(t))$ by an element of $\widetilde\Theta_0(\mathcal{H},dF(t))$.
Choose an arbitrary element $[f(t)]\in \widetilde\Theta(\mathcal{H},dF(t))$. 
Let $f(t)$ have the form~(\ref{f2_9}) ($t\in\mathbb{R}$).
Set
\begin{equation}
\label{f2_35}
f_n(t) = \sum_{k=0}^r \varphi_k^{[n]} (t) h_k,\qquad  t\in\mathbb{R},\quad n\in\mathbb{N},
\end{equation}
where $\varphi_k^{[n]}$ are constructed by the rule~(\ref{f2_31}).
Observe that $[f_n(t)] \in \widetilde\Theta_0$.
Then %
$$ \| [f] - [f_n] \|_{L_2(\mathcal{H},dF(t))}^2 = \Psi\left(
\sum_{k=0}^r (\varphi_k(t) - \varphi_k^{[n]}(t)) h_k,
\sum_{j=0}^r (\varphi_j(t) - \varphi_j^{[n]}(t)) h_j
\right) = $$
$$ = \sum_{k,j=0}^r \Psi\left(
(\varphi_k(t) - \varphi_k^{[n]}(t)) h_k,
(\varphi_j(t) - \varphi_j^{[n]}(t)) h_j
\right) \leq $$
$$ \leq
\sum_{k,j=0}^r \left|\Psi\left(
(\varphi_k(t) - \varphi_k^{[n]}(t)) h_k,
(\varphi_k(t) - \varphi_k^{[n]}(t)) h_k
\right)\right|^{\frac{1}{2}} * $$
$$ * 
\left|\Psi\left(
(\varphi_j(t) - \varphi_j^{[n]}(t)) h_j,
(\varphi_j(t) - \varphi_j^{[n]}(t)) h_j
\right)\right|^{\frac{1}{2}} 
=
$$
$$ = 
\sum_{k,j=0}^r \left|
\int_{\mathbb{R}}
|\varphi_k(t) - \varphi_k^{[n]}(t)|^2 d(F(t) h_k,h_k)_{\mathcal{H}}
\right|^{\frac{1}{2}} * $$
$$ * \left|
\int_{\mathbb{R}}
|\varphi_j(t) - \varphi_j^{[n]}(t)|^2 d(F(t) h_j,h_j)_{\mathcal{H}}
\right|^{\frac{1}{2}} \rightarrow 0, $$
as $n\rightarrow\infty$.
Thus, $\widetilde\Theta_0(\mathcal{H},dF(t))$ is dense in $L_2(\mathcal{H},dF(t))$.

Set
$\Phi_1 = \Phi_1(\mathcal{H},dF(t)) = \{ \varphi(t)\in \Phi(\mathcal{H},dF(t)):\ t\varphi(t)\in \Phi(\mathcal{H},dF(t)) \}$.
Observe that $\Phi_0\subseteq\Phi_1$.
Denote by $\Theta_1 = \Theta_1(\mathcal{H},dF(t))$ a set of all $\mathcal{H}$-valued functions $f(t)$, $t\in\mathbb{R}$, which
admit a representation of the following form:
\begin{equation}
\label{f2_37}
f(t) = \sum_{k=0}^r \varphi_k(t) h_k,\qquad  \varphi_k\in\Phi_1,\ h_k\in\mathcal{H},\ r\in\mathbb{Z}_+.
\end{equation}
Notice that $\Theta_0\subseteq\Theta_1\subseteq\Theta$. Set
$$ \widetilde\Theta_1 = \widetilde\Theta_1(\mathcal{H},dF(t)) = 
\left\{
[f(t)],\ f(t)\in \Theta_1(\mathcal{H},dF(t))
\right\}, $$
where  $[f]$ means the class of the equivalence in $L_2(\mathcal{H},dF(t))$ which contains $f$.
Since $\widetilde\Theta_0\subseteq \widetilde\Theta_1 \subseteq \widetilde\Theta$, we conclude that
$\widetilde\Theta_1(\mathcal{H},dF(t))$ is dense in $L_2(\mathcal{H},dF(t))$.

Consider the following operator $\mathcal{A}_0$ in $L_2(\mathcal{H},dF(t))$:
\begin{equation}
\label{f2_39}
\mathcal{A}_0 \left[
f(t)
\right]
=
\left[
t f(t)
\right],\qquad  f\in\Theta_1,
\end{equation}
with the domain $D(\mathcal{A}_0) = \widetilde\Theta_1$.
Let us check that the operator $\mathcal{A}_0$ is well-defined.
Suppose that $f_1,f_2\in\Theta_1$ belong to the same class of the equivalence.
We need to verify that $[t f_1(t)] = [t f_2(t)]$. Choose an arbitrary element $g\in\Theta_1$.
By the definition of $\Psi$ it is seen that
\begin{equation}
\label{f2_40}
\Psi(t f_j(t),g(t)) = \Psi(f_j(t), t g(t)),\quad j=1,2. 
\end{equation}
Therefore
$$ \left(
[tf_1(t)] - [tf_2(t)], [g(t)]
\right)_{ L_2(\mathcal{H},dF(t)) } = \Psi(tf_1(t),g(t)) - \Psi(tf_2(t),g(t)) = $$
\begin{equation}
\label{f2_50}
= \Psi(f_1(t),tg(t)) - \Psi(f_2(t),tg(t)) =
\left(
[f_1(t)] - [f_2(t)], [tg(t)]
\right)_{ L_2(\mathcal{H},dF(t)) } = 0. 
\end{equation}
Since $\widetilde\Theta_1(\mathcal{H},dF(t))$ is dense in $L_2(\mathcal{H},dF(t))$, we may choose
a sequence $[g_n]\in \widetilde\Theta_1(\mathcal{H},dF(t))$, $n\in\mathbb{N}$, which tends to 
$[tf_1(t)] - [tf_2(t)]$, as $n\rightarrow\infty$.
Writing relation~(\ref{f2_50}) with $g=g_n$ and passing to the limit as $n\rightarrow\infty$ we
get $[t f_1(t)] = [t f_2(t)]$. Thus, $\mathcal{A}$ is well-defined.
Moreover, relation~(\ref{f2_40}) shows that $\mathcal{A}$ is a densely defined symmetric operator.

Let $f(t)$ be an arbitrary element of $\Theta_0(\mathcal{H},dF(t))$ with representation~(\ref{f2_33}).
Observe that
$$ \frac{1}{t\pm i} \varphi_k(t) \in \Phi_0,\qquad 0\leq k\leq r. $$
Therefore
$\frac{1}{t\pm i} f(t) \in \Theta_0(\mathcal{H},dF(t))$.
Consequently, we have
$$ (\mathcal{A}_0 \pm i E) \left[
\frac{1}{t\pm i} f(t)
\right]
=
[f(t)]. $$
Since $\widetilde\Theta_0(\mathcal{H},dF(t))$ is dense in $L_2(\mathcal{H},dF(t))$, we conclude that
the operator $\mathcal{A}_0$ is essentially self-adjoint.
We shall denote by $\mathcal{A} = \overline{\mathcal{A}_0}$ the corresponding self-adjoint operator in $L_2(\mathcal{H},dF(t))$.

\section{A description of solutions of the moment problem.}

First assume that the moment problem~(\ref{f1_1}) is given and it has a solution $F(t)$. Choose an arbitrary $r\in\mathbb{Z}_+$,
arbitrary elements $\{ h_k \}_{0}^r$ in $\mathcal{H}$ and write
$$ \sum_{k,l=0}^r (S_{k+l} h_k, h_l)_{\mathcal{H}} =
\sum_{k,l=0}^r \left( \int_{\mathbb{R}} t^{k+l} dF(t) h_k, h_l \right)_{\mathcal{H}} = $$
$$ =
\sum_{k,l=0}^r 
\lim\limits_{ {a\to -\infty}\atop{b\to +\infty}  } 
\left(
\int_a^b t^{k+l} dF(t) h_k, h_l \right)_{\mathcal{H}} = $$
$$ =
\sum_{k,l=0}^r 
\lim\limits_{ {a\to -\infty}\atop{b\to +\infty}  } 
\lim_{\delta\rightarrow +0}
\left(
\sum_{j=0}^{N-1} (t_j^*)^{k+l} (F(t_{j+1})-F(t_j)) h_k,
, h_l \right)_{\mathcal{H}} = $$
$$ =
\lim\limits_{ {a\to -\infty}\atop{b\to +\infty}  } 
\lim_{\delta\rightarrow +0}
\sum_{j=0}^{N-1}
\left(
(F(t_{j+1})-F(t_j)) \sum_{k=0}^r (t_j^*)^{k} h_k,
, \sum_{l=0}^r (t_j^*)^{l} h_l \right)_{\mathcal{H}} \geq 0. $$
Here $\delta$ is the diameter of a partition~(\ref{f2_16}) of $[a,b]$
and $t_j^*\in [t_j, t_{j+1})$ are arbitrary points.
Therefore condition~(\ref{f1_2}) holds.

On the other hand, suppose now that the moment problem~(\ref{f1_1}) is given and condition~(\ref{f1_2}) is satisfied.
Following the idea in~\cite{citNew_4290_SK} we may consider abstract symbols $\varepsilon_j$, $j\in\mathbb{Z}_+$, and form a formal sum $h$:
\begin{equation}
\label{f3_1}
h = \sum_{j=0}^\infty h_j \varepsilon_j, 
\end{equation}
where $h_j\in \mathcal{H}$, and all but finite number of $h_j$ are zero elements of $\mathcal{H}$.
In what follows, we shall also often write finite sums, where missing summands $h_j \varepsilon_j$ are assumed to have $h_j=0$.

For $\alpha\in \mathbb{C}$ we set 
$\alpha h = \sum_{j=0}^\infty (\alpha h_j) \varepsilon_j$.
If
\begin{equation}
\label{f3_2}
g = \sum_{j=0}^\infty g_j \varepsilon_j,\qquad g_j\in \mathcal{H},
\end{equation}
where all but finite number of $g_j$ are zero elements of $\mathcal{H}$,
then we set
$h+g =  \sum_{j=0}^\infty (h_j + g_j) \varepsilon_j$. 
Thus, the set of all formal sums of type~(\ref{f3_1}) forms a complex linear vector space which we denote by $\mathfrak{B}$.
For $h,g\in \mathfrak{B}$ with representations~(\ref{f3_1}), (\ref{f3_2}) we set
\begin{equation}
\label{f3_3}
\Phi(h,g) = \sum_{j,k=0}^\infty  ( S_{j+k} h_j, g_k )_{\mathcal{H}}. 
\end{equation}%
The functional $\Phi$ is sesquilinear and $\overline{\Phi(h,g)} = \Phi(g,h)$, $\Phi(h,h)\geq 0$.
As usual, if $\Phi(h-g,h-g)=0$, we put elements $h$ and $g$ to the same equivalence class, denoted by $[h]$ or $[g]$.
The set of all equivalence classes we denote by $\mathfrak{L}$. By the completion of $\mathfrak{L}$ we get
a Hilbert space $H$.
Set
\begin{equation}
\label{f3_4}
x_{h,j}  := [ h \varepsilon_j ],\qquad h\in\mathcal{H},\ j\in\mathbb{Z}_+.
\end{equation}
Notice that
\begin{equation}
\label{f3_5}
( x_{h,j}, x_{g,k} )_H = ( S_{j+k} h, g )_{\mathcal{H}},\qquad h,g\in\mathcal{H},\ j,k\in\mathbb{Z}_+.
\end{equation}
Denote 
$$ D_0 = 
\left\{ 
\sum_{k=0}^r x_{h_k,k}:\quad h_k\in \mathcal{H},\ 0\leq k\leq r;\ r\in\mathbb{Z}_+ 
\right\}. $$
Observe that $D_0$ is a linear manifold which is dense in $H$.
Consider the following linear operator $A_0$ with $D(A_0) = D_0$:
\begin{equation}
\label{f3_6}
A_0 \sum_{k=0}^r x_{h_k,k} = \sum_{k=0}^r x_{h_k,k+1},\qquad h_k\in \mathcal{H},\ r\in\mathbb{Z}_+.
\end{equation}
Let us check that $A_0$ is well-defined.
If an element $h\in D_0$ has two representations:
$$ h = \sum_{k=0}^r x_{h_k,k} = \sum_{k=0}^s x_{g_k,k},\qquad h_k,g_k\in \mathcal{H},\ r,s\in\mathbb{Z}_+, $$
then we set $\rho = \max(r,s)$, and write:
$$ \left( 
\sum_{k=0}^\rho x_{h_k,k+1} - \sum_{k=0}^\rho x_{g_k,k+1}, x_{u,l}
\right)_H
=
\left( \sum_{k=0}^{\rho} x_{h_k-g_k,k+1}, x_{u,l} \right)_H = $$
$$ = \sum_{k=0}^{\rho} ( S_{k+1+l} (h_k-g_k), u )_{\mathcal{H}} =
\left( \sum_{k=0}^{\rho} x_{h_k-g_k,k}, x_{u,l+1} \right)_H = $$
$$ = \left( \sum_{k=0}^{r} x_{h_k,k} - \sum_{k=0}^{s} x_{g_k,k}, x_{u,l+1} \right)_H 
= 0,\quad u\in\mathcal{H},\ l\in\mathbb{Z}_+, $$
where firstly appeared $h_k,g_k$ are supposed to be zero elements.
Since $\overline{D_0} = H$, we get
$\sum_{k=0}^{r} x_{h_k,k+1} = \sum_{k=0}^{s} x_{g_k,k+1}$.
Consequently, $A_0$ is well-defined.
By~(\ref{f3_5}) it is readily checked that $A_0$ is linear and symmetric. 
Notice that
\begin{equation}
\label{f3_6_1}
A_0 x_{h,j} = x_{h,j+1},\qquad h\in\mathcal{H},\ j\in\mathbb{Z}_+.
\end{equation}
Set $A = \overline{A_0}$. 
By the induction argument we get
$$ x_{h,j} = A^j x_{h,0},\qquad j\in\mathbb{Z}_+,\ h\in\mathcal{H}. $$
Let $\widetilde A\supseteq A$ be a self-adjoint operator in a Hilbert space
$\widetilde H\supseteq H$, and $\{ \widetilde E_t \}_{t\in \mathbb{R}}$ be its (strongly) left-continuous orthogonal resolution of the identity. 
We may write
$$ (S_j f,g)_{\mathcal{H}} = (x_{f,j}, x_{g,0})_H = (A^j x_{f,0},x_{g,0})_H = (\widetilde A^j x_{f,0},x_{g,0})_{\widetilde H} = $$
\begin{equation}
\label{f3_7}
j\in\mathbb{Z}_+,\ f,g\in \mathcal{H}. 
\end{equation}
Consider the following operator $I$: $\mathcal{H}\rightarrow H$:
\begin{equation}
\label{f3_8}
I h = x_{h,0},\qquad h\in \mathcal{H}.
\end{equation}
It is readily checked that $I$ is linear, and since
$$ \| Ih \|_H^2 = (x_{h,0}, x_{h,0})_H = (S_0 h,h)_{\mathcal{H}} \leq \| S_0 \| \| h \|_{\mathcal{H}}^2, $$
then $I$ is bounded.
By~(\ref{f3_7}) we obtain that
$$ (S_j f,g)_{\mathcal{H}} =
(I^* P^{\widetilde H}_H \widetilde A^j I f,g)_{\mathcal{H}},\qquad
j\in\mathbb{Z}_+,\ f,g\in \mathcal{H}. $$
Therefore
$$ S_j f =
I^* P^{\widetilde H}_H \widetilde A^j I f
=
I^* P^{\widetilde H}_H 
\lim\limits_{ {a\to -\infty}\atop{b\to +\infty}  } \int_a^b t^j d\widetilde E_t 
I f = $$
$$ =
\lim\limits_{ {a\to -\infty}\atop{b\to +\infty}  } 
I^* P^{\widetilde H}_H
\int_a^b t^j d\widetilde E_t 
I f
=
\lim\limits_{ {a\to -\infty}\atop{b\to +\infty}  } 
\int_a^b t^j d(I^* \mathbf{E}_t I) f, $$
\begin{equation}
\label{f3_9}
f\in \mathcal{H},\ j\in\mathbb{Z}_+, 
\end{equation}
where $\mathbf{E}_t$ is a strongly left-continuous spectral function of $A$ (corresponding to $\widetilde A$).
By~(\ref{f3_9}) and the properties of the orthogonal resolution of the identity it is readily checked that
the following $[\mathcal{H}]$-valued function:
\begin{equation}
\label{f3_10}
F(t) = I^* \mathbf{E}_t I - I^* \mathbf{E}_0 I,\qquad t\in \mathbb{R},
\end{equation}
is a solution of the moment problem~(\ref{f1_1}).
Thus, \textit{each strongly left-continuous spectral function of $A$ generates a solution of
the moment problem~(\ref{f1_1}) by relation~(\ref{f3_10})}.
In the meantime we proved that \textit{conditions~(\ref{f1_2}) are necessary and sufficient for the
solvability of the moment problem~(\ref{f1_1})}.

Let $\widehat F(t)$ be an arbitrary solution of the moment problem~(\ref{f1_1}). Let us check that
$\widehat F(t)$ can be constructed by relation~(\ref{f3_10}) by virtue of some strongly left-continuous spectral function $\mathbf{E}_t$ of $A$.
Consider the Hilbert space
$L_2 = L_2(\mathcal{H}, d\widehat F(t))$ and the operator $\mathcal{A}$ in it, both defined as in the previous section.
Consider the following vector polynomials:
\begin{equation}
\label{f3_12}
p(t) = \sum_{k=0}^r t^k h_k,\quad q(t) = \sum_{l=0}^q t^l g_l,\qquad h_k,g_l\in \mathcal{H},\ r,q\in\mathbb{Z}_+.
\end{equation}
Observe that these vector polynomials belong to $\Theta_1(\mathcal{H},d\widehat F(t))$.
Thus, the equivalence classes in $L_2(\mathcal{H},d\widehat F(t))$, which contain these vector polynomials, belong to 
$\widetilde\Theta_1(\mathcal{H},d\widehat F(t))$ (which is dense in $L_2(\mathcal{H},d\widehat F(t))$).
As usual in such situations, we shall say that $p,q$ belong to $L_2(\mathcal{H}, d\widehat F(t))$.
By~(\ref{f2_25}),(\ref{f3_5}) we may write
$$ ([p],[q])_{L_2(\mathcal{H}, d\widehat F(t))}
= \int_{\mathbb{R}} \left(
d\widehat F(t) p(t), q(t)
\right)_{\mathcal{H}}
= $$
$$ = \sum_{k=0}^r \sum_{l=0}^q \int_{\mathbb{R}} t^{k+l} d \left( \widehat F(t) h_k, g_l \right)_{\mathcal{H}} = 
\sum_{k=0}^r \sum_{l=0}^q ( S_{k+l} h_k, g_l )_{\mathcal{H}} = $$
\begin{equation}
\label{f3_13}
= \sum_{k=0}^r \sum_{l=0}^q (x_{h_k,k}, x_{g_l,l})_H = \left( \sum_{k=0}^r x_{h_k,k}, \sum_{l=0}^q x_{g_l,l} \right)_H. 
\end{equation} %
Denote by $P (\mathcal{H}, d\widehat F(t))$ a set of all (classes of the equivalence which contain) polynomials of
type~(\ref{f3_12})
from $L_2(\mathcal{H}, d\widehat F(t))$.
Set $L_{2;0}(\mathcal{H}, d\widehat F(t)) = \overline{P (\mathcal{H}, d\widehat F(t))}$.
Consider the following transformation:
$$ W_0 \left[
\sum_{k=0}^r t^k h_k 
\right]
= \sum_{k=0}^r x_{h_k,k},\qquad h_k\in \mathcal{H},\ r\in\mathbb{Z}_+. $$
This transformation maps $P (\mathcal{H}, d\widehat F(t))$ on the whole $D_0$ (which is dense in $H$).
Let us check that $W_0$ is well-defined. Suppose that $p,q$ from~(\ref{f3_12})
belong to the same class of the equivalence. Without loss of generality we may assume that $r=q$, since we may add to $p$ and $q$ 
any finite number of zero addends (and this addition does not change the values of $W_0 [p]$ and $W_0 [q]$). 
We may write %
$$ 0 = 
\left\|
\sum_{k=0}^r t^k h_k - \sum_{k=0}^r t^k g_k 
\right\|_{L_2(\mathcal{H}, d\widehat F(t))}^2
=
\left\|
\sum_{k=0}^r t^k (h_k-g_k) 
\right\|_{L_2(\mathcal{H}, d\widehat F(t))}^2
=
$$
$$ =
\left(
\sum_{k=0}^r t^k (h_k-g_k),
\sum_{l=0}^r t^l (h_l-g_l) 
\right)_{L_2(\mathcal{H}, d\widehat F(t))} = $$
$$ = \left( \sum_{k=0}^r x_{h_k-g_k,k}, \sum_{l=0}^r x_{h_l-g_l,l} \right)_H 
= \left\|
\sum_{k=0}^r x_{h_k,k} - \sum_{k=0}^r x_{g_k,k}
\right\|^2_H,
$$ 
and therefore
$\sum_{k=0}^r x_{h_k,k} = \sum_{k=0}^r x_{g_k,k}$.
Thus, $W_0$ is well-defined.
Moreover, $W_0$ is linear and relation~(\ref{f3_13}) shows that $W_0$ is isometric. By continuity we extend $W_0$ to an isometric
transformation $W$ which maps $L_{2;0}(\mathcal{H}, d\widehat F(t))$ on the whole $H$.
Observe that
\begin{equation}
\label{f3_14}
A_0 x_{h,j} = x_{h,j+1} = W \left[ t^{j+1} h \right] = W\mathcal{A} \left[ t^j h \right] = W \mathcal{A} W^{-1} x_{h,j}, 
\end{equation}
for arbitrary $h\in \mathcal{H}$, $j\in\mathbb{Z}_+$. 
Set
$$ L_{2;1}(\mathcal{H}, d\widehat F(t)) = 
L_2(\mathcal{H}, d\widehat F(t))\ominus L_{2;0}(\mathcal{H}, d\widehat F(t)), $$
and $\mathbf W = W\oplus E_{L_{2;1}(\mathcal{H}, d\widehat F(t))}$.
Observe that $\mathbf W$ is an isometric transformation which maps $L_2(\mathcal{H}, d\widehat F(t))$
on $H_1 := H\oplus L_{2;1}(\mathcal{H}, d\widehat F(t))$.
Set $\widetilde A := \mathbf W \mathcal{A} \mathbf{W}^{-1}$. By~(\ref{f3_14}) we conclude that $\widetilde A$ is a self-adjoint
extension of $A_0$.
Let $\{ \widetilde E_t \}_{t\in \mathbb{R}}$ be the strongly left-continuous resolution of the identity of $\widetilde A$.
Denote by $\mathbf E_t$ ($t\in \mathbb{R}$) the corresponding strongly left-continuous spectral function of $A$.
For arbitrary $h,g\in \mathcal{H}$, $z\in\mathbb{R}_e$ we may write %
$$ \int_{\mathbb{R}} \frac{1}{t-z} d(I^* \mathbf{E}_t I h,g)_{\mathcal H} = 
\int_{\mathbb{R}} \frac{1}{t-z} d(\widetilde E_t Ih, Ig)_{H_1} = $$
$$ =
\left( 
\left( \widetilde A -z E_{H_1}
\right)^{-1} \mathbf W [h], \mathbf W [g] 
\right)_{H_1} = 
\left( 
\mathbf W^{-1} \left( \widetilde A - z E_{H_1}
\right)^{-1} \mathbf W [h], [g] 
\right)_{L_2} =
$$
$$ = \left( 
\left( \mathcal{A} - z E_{L_2} 
\right)^{-1} [h], [g] 
\right)_{L_2} =
\left( 
\left[ 
\frac{1}{t-z} h
\right], [g] 
\right)_{L_2} = $$
$$ = \int_{\mathbb{R}} 
\frac{1}{t-z} d(\widehat F(t)h, g)_{\mathcal{H}}. $$
By the Stieltjes-Perron inversion formula we conclude that $\widehat F(t) = I^* \mathbf E_t I - I^* \mathbf{E}_0 I$.
\begin{thm}
\label{t3_1}
Let the operator Hamburger moment problem~(\ref{f1_1}) be given and
condition~(\ref{f1_2}) holds. Let an operator $A_0$ in a Hilbert space $H$ be constructed as in~(\ref{f3_6}), 
$A = \overline{A_0}$.
All solutions of the moment problem have the following form:
\begin{equation}
\label{f3_15}
F(t) = I^* \mathbf{E}_t I - I^* \mathbf{E}_0 I,\qquad t\in\mathbb{R},
\end{equation}
where $I$ is defined by~(\ref{f3_8}) and $\mathbf{E}_t$ is a strongly left-continuous spectral function of $A$.
On the other hand, each strongly left-continuous spectral function of $A$ generates by~(\ref{f3_15})
a solution of the moment problem.
Moreover, the correspondence between all strongly left-continuous spectral functions of $A$ and
all solutions of the moment problem~(\ref{f1_1}) is one-to-one.
\end{thm}
\textbf{Proof.}
It remains to verify that different strongly left-continuous spectral functions of $A$ generate
different solutions of the moment problem~(\ref{f1_1}).
Set
$L_0 := \{ x_{h,0} \}_{h\in \mathcal{H}}$.
Choose an arbitrary element $x\in\mathfrak{L}$. Observe that $x$ has the following form:
$x = \sum_{j=0}^r x_{h_j,j}$, $h_j\in\mathcal{H}$, with some $r\in\mathbb{N}$.
For arbitrary $z\in\mathbb{R}_e$ there exists the following representation:
\begin{equation}
\label{f3_16}
x = v + y,\qquad v\in \mathcal{M}_z(A),\ y\in L_0.
\end{equation}
Here $v$ and $y$ may depend on the choice of $z$.
In fact, choose an element $u\in D(A_0)$ of the following form:
$u = \sum_{k=0}^{r-1} x_{g_k,k}$, with arbitrarily chosen $g_k\in\mathcal{H}$.
We may write:    %
$$ (A  - z E_H) u = A_0 u - z u = \sum_{k=0}^{r-1} x_{g_k,k+1} - \sum_{k=0}^{r-1} z x_{g_k,k} = 
\sum_{k=1}^r x_{g_{k-1},k} -  \sum_{k=0}^{r-1} z x_{g_k,k} = $$
\begin{equation}
\label{f3_17}
= \sum_{k=1}^{r-1} x_{g_{k-1} - z g_k,k} + x_{g_{r-1},r} + (-z) x_{g_0,0}, 
\end{equation}
where the sum on the right is empty if $r=1$.
Set
$$ g_{r-1} = h_r, $$
\begin{equation}
\label{f3_19}
g_{k-1} = z g_k + h_k,\qquad  k = r-1,r-2,...,1.
\end{equation}
In the case $r=1$ equalities~(\ref{f3_19}) are redundant.
Thus, we have defined $g_{r-1},g_{r-2},...,g_0$.

Consider $u$ with this choice of $g_k$ and set $v = (A - z E_H) u\in \mathcal{M}_z(A)$.
By~(\ref{f3_17}),(\ref{f3_19}) we see that
$$ x = v + (-1) x_{(-1)h_0 + (-z)g_0, 0}. $$
Therefore relation~(\ref{f3_16}) is proved.     %

Suppose to the contrary that two different strongly left-continuous spectral functions
$\mathbf E_{j,t}$, $j=1,2$, generate the same solution of the moment problem:
\begin{equation}
\label{f3_20}
I^* \mathbf E_{1,t} I - I^* \mathbf E_{1,0} I = I^* \mathbf E_{2,t} I - I^* \mathbf E_{2,0} I,\qquad t\in\mathbb{R}.
\end{equation}
For arbitrary $f\in \mathcal{H}$ and $t\in\mathbb{R}$ we may write:
$$ (\mathbf E_{1,t} x_{f,0}, x_{f,0})_H - (\mathbf E_{1,0} x_{f,0}, x_{f,0})_H = (\mathbf E_{2,t} x_{f,0}, x_{f,0})_H -
(\mathbf E_{2,0} x_{f,0}, x_{f,0})_H. $$
Multiplying by $\frac{1}{t-z}$ and integrating we get
$$ \int_{\mathbb{R}} \frac{1}{t-z} d(\mathbf E_{1,t} If, If)_H = 
\int_{\mathbb{R}} \frac{1}{t-z} d(\mathbf E_{2,t} If, If)_H,\qquad z\in\mathbb{R}_e. $$
Therefore
$I^* \mathbf R_{1,z} I = I^* \mathbf R_{2,z} I$,
where $\mathbf R_{j,z}$ is the generalized resolvent corresponding to $\mathbf E_{j,t}$, $j=1,2$.
Then
\begin{equation}
\label{f3_22}
(\mathbf R_{1,z} x_{f,0}, x_{g,0})_H = (\mathbf R_{2,z} x_{f,0}, x_{g,0})_H,\qquad f,g\in\mathcal{H},\ z\in\mathbb{R}_e.
\end{equation}
Let $\mathbf R_{j,z}$ be generated by a self-adjoint extension $\widetilde A_j$ of $A$ in a Hilbert space
$\widetilde H_j\supseteq H$, $j=1,2$. 
Since for arbitrary $f\in D(A), z\in\mathbb{R}_e$ and $j=1,2$ we have
$\left( \widetilde A_j - z E_{\widetilde H_j} \right)^{-1} \left(A - z E_H \right) f =
f\in H$, 
then
\begin{equation}
\label{f3_24}
\mathbf R_{1,z} u = \mathbf R_{2,z} u,\qquad u\in \mathcal{M}_z(A),\ z\in\mathbb{R}_e.
\end{equation}
Observe that
$$ (\mathbf R_{j,z} x_{f,0}, u) = (x_{f,0}, \mathbf R_{j,z}^* u) = 
\left( x_{f,0}, \mathbf R_{ j, \overline{z} } u \right), $$
where $f\in\mathcal{H}$, $u\in \mathcal{M}_{ \overline{z} }(A)$, $j=1,2$; $z\in\mathbb{R}_e$. 
By~(\ref{f3_24}) we get
\begin{equation}
\label{f3_26}
(\mathbf R_{1,z} x_{f,0}, u)  = (\mathbf R_{2,z} x_{f,0}, u),\qquad u\in \mathcal{M}_{ \overline{z} }(A),\ f\in \mathcal{H},\ 
z\in\mathbb{R}_e.
\end{equation}
Choose an arbitrary element $w\in\mathfrak{L}$ and $z\in \mathbb{R}_e$. By~(\ref{f3_16}) we may write:
$w = v + y$, where $v\in \mathcal{M}_{ \overline{z} }(A)$, $y\in L_0$. 
By~(\ref{f3_22}),(\ref{f3_26}) we obtain that
$(\mathbf R_{1,z} x_{f,0}, w)  = (\mathbf R_{2,z} x_{f,0}, w)$, $\forall f\in\mathcal{H}$. Therefore
\begin{equation}
\label{f3_28}
\mathbf R_{1,z} x = \mathbf R_{2,z} x,\qquad x\in L_0,\ z\in\mathbb{R}_e.
\end{equation}
For an arbitrary $w\in\mathfrak{L}$ using~(\ref{f3_16}) we may write:
$w = v' + y'$, where $v'\in \mathcal{M}_z(A)$, $y'\in L_0$, $z\in\mathbb{R}_e$.
By~(\ref{f3_24}),(\ref{f3_28}) we get
$\mathbf R_{1,z} w = \mathbf R_{2,z} w$, $z\in\mathbb{R}_e$.
Therefore $\mathbf E_{1,t}=\mathbf E_{2,t}$. This contradiction completes the proof.
$\Box$ %

By~(\ref{f3_15}) we see that each solution $F(t)$ of the moment problem satisfies the following relation:
\begin{equation}
\label{f3_30}
(F(t)h,h)_{\mathcal{H}} = (\mathbf{E}_t Ih, Ih) - (\mathbf{E}_0 Ih, Ih),\qquad h\in\mathcal{H},
\end{equation} 
where $\mathbf{E}_t$ is a strongly left-continuous spectral function of the operator $A$.
From this relation it follows that
\begin{equation}
\label{f3_32}
\int_{\mathbb{R}} \frac{1}{t-z} d(F(t)h,h)_{\mathcal{H}} = 
(\mathbf{R}_z x_{h,0}, x_{h,0})_H,\qquad h\in\mathcal{H},\ 
z\in \mathbb{R}_e,
\end{equation}
where $\mathbf{R}_z$ is the generalized resolvent of $A$ which corresponds to $\mathbf{E}_t$.

We need to shift the vector $x_{h,0}$ suitably, in order to fit into the block 
of the operator Frobenius formula in the further construction.
Set
$$ y_{h,0} := (A - iE_H) x_{h,0} = x_{h,1} - i x_{h,0},\qquad h\in\mathcal{H}. $$
Let $\widetilde A$ be a self-adjoint operator in a Hilbert space $\widetilde{H}\supseteq H$ which generates
the generalized resolvent $\mathbf{R}_z$.
Notice that
$$ (\mathbf{R}_z x_{h,0}, x_{h,0})_H =
(R_z(\widetilde A) R_i(\widetilde A) y_{h,0}, x_{h,0})_H = $$
$$ = 
\frac{1}{z-i}
\left(
(R_z(\widetilde A) y_{h,0}, x_{h,0})_H -(R_i(\widetilde A) y_{h,0}, x_{h,0})_H 
\right) =
$$
\begin{equation}
\label{f3_200}
= \frac{1}{z-i} (R_z(\widetilde A) y_{h,0}, x_{h,0})_H
- \frac{1}{z-i} (S_0 h,h)_{\mathcal{H}},\qquad h\in\mathcal{H},\ z\in \mathbb{R}_e\backslash\{ i \}.
\end{equation}
We may write:
$$ (R_z(\widetilde A) y_{h,0}, x_{h,0})_H =
(R_z(\widetilde A) y_{h,0}, R_i(\widetilde A) y_{h,0})_H = 
(R_{-i}(\widetilde A) R_z(\widetilde A) y_{h,0}, y_{h,0})_H = $$
$$ = 
-\frac{1}{z+i}
\left(
(R_{-i}(\widetilde A) y_{h,0}, y_{h,0})_H - (R_z(\widetilde A) y_{h,0}, y_{h,0})_H
\right) =
$$
$$ =
\frac{1}{z+i}
(R_z(\widetilde A) y_{h,0}, y_{h,0})_H
-
\frac{1}{z+i}
( y_{h,0}, x_{h,0} )_H
= $$
\begin{equation}
\label{f3_210}
=
\frac{1}{z+i}
(R_z(\widetilde A) y_{h,0}, y_{h,0})_H
-
\frac{1}{z+i}
( (S_1 - i S_0) h,h)_{\mathcal{H}},\qquad h\in\mathcal{H},\ z\in \mathbb{R}_e\backslash\{ -i \}.
\end{equation}
By~(\ref{f3_200}),(\ref{f3_210}) we get
$$ (\mathbf{R}_z x_{h,0}, x_{h,0})_H =
\frac{1}{z^2 + 1} (\mathbf{R}_z y_{h,0}, y_{h,0})_H -
\frac{1}{z^2 + 1} 
( (z S_0 + S_1) h,h)_{\mathcal{H}},\ $$
$$ h\in\mathcal{H},\ z\in \mathbb{R}_e\backslash\{ -i, i \}. $$
Relation~(\ref{f3_32}) takes the following form:                       %
$$  \int_{\mathbb{R}} \frac{1}{t-z} d(F(t)h,h)_{\mathcal{H}} =  $$
\begin{equation}
\label{f3_230}
= \frac{1}{z^2 + 1} (\mathbf{R}_z y_{h,0}, y_{h,0})_H -
\frac{1}{z^2 + 1} 
( (z S_0 + S_1) h,h)_{\mathcal{H}},\qquad h\in\mathcal{H},\ 
z\in \mathbb{R}_e\backslash\{ -i, i \}.
\end{equation}
Consider Cayley's transformation of the closed symmetric operator $A$:
$V := (A+iE_H)(A-iE_H)^{-1}$.
Generalized resolvents of $A$ and $V$ are related in the following way:
\begin{equation}
\label{f3_240}
(1-\zeta) \mathbf{R}_\zeta (V) = E_H + (z-i) \mathbf{R}_z(A),\qquad \zeta = \frac{z-i}{z+i}(\in\mathbb{D}),\  z\in\mathbb{C}_+.
\end{equation}
This relation follows from the fact that the corresponding resolvents of $\widetilde A$ and
$\widetilde V := (\widetilde A+iE_H)(\widetilde A-iE_H)^{-1}$ are related in the same manner
(and by applying $P^{\widetilde H}_H$).
This fact was noticed in~\cite[pp. 370-371]{citMFAT_2200_Ch}.   %

Extracting $\mathbf{R}_z(A)$ from~(\ref{f3_240}), substituting into relation~(\ref{f3_230}) and
simplifying we get
$$  \int_{\mathbb{R}} \frac{1}{t-z} d(F(t)h,h)_{\mathcal{H}} =  $$
$$ = \frac{2i}{(z^2 + 1)^2} (\mathbf{R}_{\frac{z-i}{z+i}}(V) y_{h,0}, y_{h,0})_H -
\frac{1}{(z-i)(z^2 + 1)}
( (S_2 + S_0) h,h)_{\mathcal{H}} - $$
\begin{equation}
\label{f3_250}
-
\frac{1}{z^2 + 1} 
( (z S_0 + S_1) h,h)_{\mathcal{H}},\qquad h\in\mathcal{H},\ 
z\in \mathbb{C}_+\backslash\{ i \}.
\end{equation}

Applying Chumakin's formula for the generalized resolvents of a closed isometric operator (\cite{cit_950_Ch}) we obtain that
each solution $F(t)$ of the moment problem~(\ref{f1_1}) satisfies the following relation:
$$  \int_{\mathbb{R}} \frac{1}{t-z} d(F(t)h,h)_{\mathcal{H}} =  $$
$$ = \frac{2i}{(z^2 + 1)^2} \left(
\left[
E_H
-
\frac{z-i}{z+i}
\left(
V\oplus
\Phi_{\frac{z-i}{z+i}}
\right)
\right]^{-1}
y_{h,0}, y_{h,0}
\right)_H - $$
$$ - \frac{1}{(z-i)(z^2 + 1)}
( (S_2 + S_0) h,h)_{\mathcal{H}} -
\frac{1}{z^2 + 1} 
( (z S_0 + S_1) h,h)_{\mathcal{H}},\qquad h\in\mathcal{H},\ $$
\begin{equation}
\label{f3_260}
 z\in \mathbb{C}_+\backslash\{ i \},
\end{equation}
where $\Phi_\zeta\in \mathcal{S}(\mathbb{D};H\ominus D(V),H\ominus R(V))$.

Denote the expression in the square brackets in~(\ref{f3_260}) by $M=M_z(\Phi)$.
Applying to $M$ the operator Frobenius inversion formula in~\cite[Proposition 1]{citNew_1700_Z} for the
decomposition $H = \mathcal{M}_i(A) \oplus \mathcal{N}_i(A)$ we get:
\begin{equation}
\label{f3_270}
M^{-1} =
\left(
\begin{array}{cc} \mathbf{A}_z^{-1} + \mathbf{A}_z^{-1} \mathbf{B}_z(\Phi) (\mathcal{H}_z(\Phi))^{-1} \mathbf{C}_z \mathbf{A}_z^{-1} & \ast \\
\ast & \ast \end{array}
\right),
\end{equation}
where
\begin{equation}
\label{f3_290}
\widehat{\mathbf{A}}_z := \mathbf{A}_z^{-1} = \left(- \frac{z+i}{z-i} \right) 
\left(
P^H_{\mathcal{M}_i(A)} V -
\frac{z+i}{z-i}
E_{\mathcal{M}_i(A)}
\right)^{-1},
\end{equation}
\begin{equation}
\label{f3_295}
\mathbf{B}_z(\Phi) = 
\left(
- \frac{z-i}{z+i} 
\right)
P^H_{\mathcal{M}_i(A)} \Phi_{\frac{z-i}{z+i}},
\end{equation}
\begin{equation}
\label{f3_310}
\mathbf{C}_z = 
\left(
- \frac{z-i}{z+i}
\right)
P^H_{\mathcal{N}_i(A)} V,\quad
\end{equation}
\begin{equation}
\label{f3_315}
\mathbf{D}_z(\Phi) = E_{\mathcal{N}_i(A)} + 
\left( - \frac{z-i}{z+i} \right)
P^H_{\mathcal{N}_i(A)} \Phi_{\frac{z-i}{z+i}},
\end{equation}
and
$$ \mathcal{H}_z(\Phi) = 
\mathbf{D}_z(\Phi) - \mathbf{C}_z \mathbf{A}_z^{-1} \mathbf{B}_z(\Phi),
$$
\begin{equation}
\label{f3_320}
z\in \mathbb{C}_+\backslash\{ i \},\ \Phi_\zeta\in \mathcal{S}(\mathbb{D};H\ominus D(V),H\ominus R(V)).
\end{equation}
By $\ast$ we denote the blocks which are not of interest for us.

\noindent
Relation~(\ref{f3_260}) takes the following form:
$$  \int_{\mathbb{R}} \frac{1}{t-z} d(F(t)h,h)_{\mathcal{H}} =  $$
$$ = \frac{2i}{(z^2 + 1)^2} \left(
\left(
\widehat{\mathbf{A}}_z + \widehat{\mathbf{A}}_z \mathbf{B}_z(\Phi) (\mathcal{H}_z(\Phi))^{-1} \mathbf{C}_z \widehat{\mathbf{A}}_z
\right)
y_{h,0}, y_{h,0}
\right)_H - $$
$$ - \frac{1}{(z-i)(z^2 + 1)}
( (S_2 + S_0) h,h)_{\mathcal{H}} -
\frac{1}{z^2 + 1} 
( (z S_0 + S_1) h,h)_{\mathcal{H}},\qquad h\in\mathcal{H},\ $$
\begin{equation}
\label{f3_260_1}
 z\in \mathbb{C}_+\backslash\{ i \},
\end{equation}
where $\Phi_\zeta\in \mathcal{S}(\mathbb{D};H\ominus D(V),H\ominus R(V))$.

Consider the following operator $K$, which maps $\mathcal{H}$ into $\mathcal{M}_i(A)$:
\begin{equation}
\label{f3_350}
K h = y_{h,0} = x_{h,1} - i x_{h,0},\qquad h\in\mathcal{H}.
\end{equation}
It is readily checked that $K$ is linear, and
$$ \| Kh \|^2 = (x_{h,1} - i x_{h,0}, x_{h,1} - i x_{h,0}) = (S_2 h,h) + (S_0 h,h) \leq 
(\| S_0 \| + \| S_2 \|) \| h \|^2, $$
for an arbitrary $h\in\mathcal{H}$.
Therefore $K$ is bounded.

Relation~(\ref{f3_260_1}) may be rewritten in the following form:
$$  \int_{\mathbb{R}} \frac{1}{t-z} d(F(t)h,h)_{\mathcal{H}} =  $$
$$ = \frac{2i}{(z^2 + 1)^2} \left(
K^*
\left(
\widehat{\mathbf{A}}_z + \widehat{\mathbf{A}}_z \mathbf{B}_z(\Phi) (\mathcal{H}_z(\Phi))^{-1} \mathbf{C}_z \widehat{\mathbf{A}}_z
\right)
K h, h
\right)_H - $$
$$ - \frac{1}{(z-i)(z^2 + 1)}
( (S_2 + S_0) h,h)_{\mathcal{H}} -
\frac{1}{z^2 + 1} 
( (z S_0 + S_1) h,h)_{\mathcal{H}},\qquad h\in\mathcal{H},\ $$
\begin{equation}
\label{f3_260_2}
 z\in \mathbb{C}_+\backslash\{ i \},
\end{equation}
where $\Phi_\zeta\in \mathcal{S}(\mathbb{D};H\ominus D(V),H\ominus R(V))$.

\begin{thm}
\label{t3_2}
Let the operator Hamburger moment problem~(\ref{f1_1}) be given and
condition~(\ref{f1_2}) holds. Let an operator $A$ in a Hilbert space $H$ be constructed as in~(\ref{f3_6}), 
an operator $K$ be constructed by~(\ref{f3_350}), $A = \overline{A_0}$, $V = (A+iE_H)(A-iE_H)^{-1}$. 
Define $\widehat{\mathbf{A}}_z$, $\mathbf{B}_z(\Phi)$, $\mathbf{C}_z$, $\mathbf{D}_z(\Phi)$,
$\mathcal{H}_z(\Phi)$ by relations~(\ref{f3_290})-(\ref{f3_320}).
Each solution $F(t)$ of the moment problem~(\ref{f1_1}) satisfies relation~(\ref{f3_260_2}) with some
function $\Phi_\zeta\in \mathcal{S}(\mathbb{D};H\ominus D(V),H\ominus R(V))$.
On the other hand, for each function $\Phi_\zeta\in \mathcal{S}(\mathbb{D};H\ominus D(V),H\ominus R(V))$
there exists a solution $F(t)$ of the moment problem~(\ref{f1_1}) which satisfies~(\ref{f3_260_2}).
\end{thm}
\textbf{Proof.}
The proof follows from the considerations before the statement of the theorem.
$\Box$

\begin{center}
{\large\bf The Nevanlinna-type parametrization for the operator Hamburger moment problem.}
\end{center}
\begin{center}
{\bf S.M. Zagorodnyuk}
\end{center}

In this paper we obtain a Nevanlinna-type parametrization for the operator Hamburger moment
problem. 
The moment problem is not supposed to be necessarily completely indeterminate.
No assumptions besides the solvability of the moment problem are posed.

}

\begin{thebibliography}{11}

\bibitem{cit_510_I}
G. M. Ilmushkin,
On solutions of the operator power moment problem,
Funktsionalnij analiz (Ulyanovsk),
\textbf{14} (1980), 74--79. (Russian). 

\bibitem{citNew_20_K}
A. Kheifets,
Hamburger moment problem: Parseval equality and A-singularity,
Journal of Functional Analysis,
\textbf{141} (1996),
374--420. 


\bibitem{cit_600_Akh}
N. I. Akhiezer,
The Classical Moment Problem and Some Related Questions in Analysis,
Oliver \& Boyd,
Edinburgh,
1965.

\bibitem{citNew_25_LN}
L. Lemnete-Ninulescu, Hamburger and Stieltjes moment problems for operators, 
ISRN Mathematical Analysis, \textbf{2014} (2014), Article ID 836839, 7 pages. 

\bibitem{cit_700_Ber}
Ju. M. Berezanskii,
Expansions in Eigenfunctions of Selfadjoint
Operators,
Amer. Math. Soc.,
Providence, RI,
1968. (Russian edition: Naukova Dumka, Kiev, 1965).

\bibitem{cit_800_AK}
N. Akhiezer, M. Krein,
Some Questions in the Theory of Moments,
Nauchno-Tehnicheskoe Izdatelstvo Ukrainy,
Kharkov,
1938. (Russian).

\bibitem{citNew_30_K}
M. G. Krein,
The fundamental propositions of the theory of representations of Hermitian operators with
deficiency index $(m,m)$,
Ukr. mat. zhurnal,
\textbf{1}, no.~2 (1949), 3--66. (Russian).



\bibitem{citNew_50_CH}
G.-N. Chen, Y.-J. Hu, The Nevanlinna-Pick interpolation problems and power moment problems for matrix-valued
functions III: The infinitely many data case,
Linear Algebra and its Applications,
\textbf{306} (2000),  59--86.

\bibitem{cit_400_L}
P. Lopez-Rodriguez,
The Nevanlinna parametrization for a matrix moment problem,
Math. Scand.,
\textbf{89} (2001),
245--267.

\bibitem{cit_970_B}
C. Berg,
The matrix moment problem,
Coimbra lecture notes,
Coimbra,
2010.


\bibitem{cit_850_Z}
S. M. Zagorodnyuk,
A description of all solutions of the matrix Hamburger moment problem in a general case,
Methods Funct. Anal. Topology,
\textbf{16}, no.~3 (2010), 271--288.

\bibitem{citNew_100_Z}
S. M. Zagorodnyuk, The Nevanlinna-type formula for the matrix Hamburger moment problem,
Methods Funct. Anal. Topology,
\textbf{18}, no.~4 (2012),  387--400.


\bibitem{citNew_300_KK}
M. G. Krein, M. A. Krasnoselskii, Fundamental theorems on the extension of Hermitian operators and certain of their applications to the theory of 
orthogonal polynomials and the problem of moments, Uspehi Matem. Nauk, \textbf{2}, no. 3(19) (1947), 60--106 (Russian).


\bibitem{citNew_500_SN}
B. Sz.-Nagy, A moment problem for self-adjoint operators, Acta Math. Acad. Sci. Hungar., \textbf{3} (1952), 285--293.

\bibitem{citNew_700_MN}
J. S. MacNerney, Hermitian moment sequences, Trans. Amer. Math. Soc., \textbf{103} (1962), 45--81.

\bibitem{citNew_900_L}
D. Leviatan,  A generalized moment problem for self-adjoint operators, Israel J. Math. \textbf{4} (1966), 113--118.



\bibitem{cit_520_IA}
G. M. Ilmushkin, E. L. Aleksandrov,
On solutions of the operator power moment problem,
Izvestiya vuzov, Matematika,
no.~3 (1989), 18--24. (Russian). 

\bibitem{cit_530_I}
G. M. Ilmushkin,
The truncated operator power moment problem,
Funktsionalnij analiz (Ulyanovsk),
\textbf{1} (1973), 51--62. (Russian). 



\bibitem{citNew_1000_A}
T. Ando, Truncated moment problems for operators,
Acta Sci. Math. (Szeged), \textbf{31}, no. 4 (1970), 319--334.

\bibitem{citNew_1400_A}
Y. M. Arlinski\u\i, Truncated Hamburger moment problem for an operator measure with compact support,
Math. Nachr., \textbf{285}, No. 14-15 (2012), 1677--1695. 

\bibitem{citNew_1700_Z}
S. Zagorodnyuk, On the truncated operator trigonometric moment problem, Concr. Oper., \textbf{2} (2015), 37--46.


\bibitem{citNew_12000_Neumark}                                      
M. A. Neumark, Spectral functions of a symmetric operator,
Izvestiya AN SSSR, \textbf{4} (1940), 277--318. (Russian).

\bibitem{citNew_14000_Neumark}                                     
M. A. Neumark, On spectral functions of a symmetric operator,
Izvestiya AN SSSR, \textbf{7} (1943), 285--296. (Russian).

\bibitem{citNew_4280_SK} 
B. Sz\"okefalvi-Nagy, A. Koranyi, Relations d'un
probl\'eme de Nevanlinna et Pick avec la th\'eorie de l'espace hilbertien,
Acta Math. Acad. Sci. Hungar., \textbf{7} (1957), 295--302.

\bibitem{citNew_4290_SK} 
B. Sz\"okefalvi-Nagy, A. Koranyi,
Operatortheoretische Behandlung und Verallgemeinerung eines  Problemkreises
in der komplexen Funktionentheorie, 
Acta Math., \textbf{100} (1958), 171--202.


\bibitem{cit_950_Ch}
M.~E. Chumakin, Generalized resolvents of isometric operators,
Sibirskiy matem. zhurnal,  {\bf 8}, no. 4  (1967), 876--892. (Russian).


\bibitem{citMFAT_2200_Ch}
M.~E. Chumakin,
On a class of the generalized resolvents of an isometric operator,
Uchen. zap. Ulyanovskogo pedinstituta,
\textbf{20}, no.~4
(1966), 361--373. (in Russian).






 
%






\end{thebibliography}
\end{document}